\documentclass[12pt,a4paper]{amsart}
\usepackage[top=35mm, bottom=30mm, left=30mm, right=30mm]{geometry}
\usepackage{mathptmx}
\usepackage{mathrsfs}
\usepackage{verbatim}
\usepackage{url}
\usepackage[all]{xy}
\usepackage{xcolor}
\usepackage[colorlinks=true,citecolor=blue]{hyperref}
\usepackage{amsmath,amsthm}
\usepackage{amssymb}
\usepackage{enumerate}
\usepackage{graphicx}
\usepackage{tikz}
\usetikzlibrary{arrows}
\usetikzlibrary{patterns}

\newtheorem{thm}{Theorem}[section]

\newtheorem{ques}[thm]{Question}

\theoremstyle{definition}

\numberwithin{equation}{section}
\begin{document}


\baselineskip=20pt


\title{Linearizations of periodic point free distal homeomorphisms on the annulus}

\author[E.~Shi]{Enhui Shi}
\address{School of Mathematics and Sciences, Soochow University, Suzhou, Jiangsu 215006, China}
\email{ehshi@suda.edu.cn}

\author[H.~Xu]{Hui Xu*}
\thanks{*Corresponding author}
\address{ Department of mathematics, Shanghai normal university, Shanghai 200234,
 China}
\email{huixu@shnu.edu.cn}

\author[Z.~Yu]{Ziqi YU}
\address{School of Mathematics and Sciences, Soochow University, Suzhou, Jiangsu 215006, China}
\email{20204207013@stu.suda.edu.cn}

\begin{abstract}
Let $\mathbb{A}$ be an annulus in the plane $\mathbb R^2$ and $g:\mathbb{A}\rightarrow \mathbb{A}$ be a boundary components preserving homeomorphism which is distal and has no periodic points.
\medskip

In \cite{SXY}, the authors show that there is a continuous decomposition $\mathcal P$ of $\mathbb{A}$ into $g$-invariant circles
 such that all the restrictions of $g$ on them share a common  irrational rotation number (also called the rotation number of $g$) and all these circles are linearly ordered by the inclusion relation on the sets of bounded components of their complements in $\mathbb R^2$.
\medskip

In this note, we show that if the decomposition  $\mathcal P$  above has a continuous section, then $g$ can be linearized, that is it
is topologically conjugate to a rigid rotation on $\mathbb{A}$. For every irrational number $\alpha\in (0, 1)$, we show the existence of such a distal homeomorphism $g$ on $\mathbb{A}$ that it cannot be linearized and its rotation number is $\alpha$.
\end{abstract}

\subjclass[2010]{}

\maketitle

\pagestyle{myheadings} \markboth{E. Shi, H. Xu, and Z. Yu}{The structure of periodic point free distal homeomorphisms}

\section{Introduction}

We refer the readers to \cite{SXY} for notions, notations, and backgrounds. In  \cite{SXY}, we obtained the following theorem.

\begin{thm}\label{main}
Let $\mathbb{A}$  be an annulus in the plane $\mathbb R^2$ and $g:\mathbb{A}\rightarrow \mathbb{A}$ be a boundary components preserving homeomorphism
which is distal and has no periodic points.
Then there is a continuous decomposition $\mathcal P$ of $\mathbb{A}$ into $g$-invariant circles
 such that all the restrictions of $g$ on them share a common  irrational rotation number and all these circles are
linearly ordered by the inclusion relation on the sets of bounded components of their complements in $\mathbb R^2$.
\end{thm}

The following question is natural.

\begin{ques}
Are the homeomorphisms $g$ in Theorem \ref{main} conjugate to rigid rotations on $\mathbb{A}$; that is, can $g$ be linearized?
\end{ques}

We deal with this question in the following results.

\begin{thm}\label{main2}
Let $\mathbb{A}$, $g$, and $\mathcal P$ be as in Theorem \ref{main}.  If there is
a continuous injection $\phi:[0,1]\rightarrow \mathbb{A}$ such that $\phi([0,1])\cap C$ has exactly one point for each $C\in\mathcal P$,
then $g$ is topologically conjugate to a rigid rotation on $\mathbb{A}$.
\end{thm}

\begin{thm}\label{example}
For each irrational number $\alpha\in (0, 1)$, there is a distal homeomorphism $g$ as in  Theorem \ref{main} which has rotation number $\alpha$
and cannot be linearized.
\end{thm}

The following question is left.

\begin{ques}\label{smooth}
Can the homeomorphisms $g$ in Theorem \ref{example} be taken to be smooth?
\end{ques}

The answer to Question \ref{smooth} may depend on the approximation properties of $\alpha$ by rational numbers.
So it more or less has some connections with Herman's question (see e.g. \cite[Chap. 4.3]{BZ}).

\section{Proof of Theorem \ref{main2}}

Let $\mathcal{P}$ be the decomposition of the annulus into minimal circles (see \cite{SXY}). A {\it transversal} of $\mathcal{P}$ is an arc in the annulus that intersects each member in $\mathcal{P}$ exactly once. Now suppose that $\gamma$ is a transversal for $\mathcal{P}$.  For each $n\in\mathbb{Z}$, let $L_{n}=g^{n}(\gamma)$.

 Recall that the map $[1,2]\rightarrow \mathcal{P}, \alpha\mapsto C_{\alpha}$  is continuous. We parametrize $\mathcal{P}$ with $C_{\alpha}: \alpha\in[1,2]$. For each $n\in\mathbb{Z}$ and $\alpha\in[1,2]$, set $\{x_{n}^{\alpha}\}=L_n\cap C_{\alpha}$. Then it is clear that $\{x_{n}^{\alpha}\}_{n\in\mathbb{Z}}$ is dense in $C_{\alpha}$. WLOG, we may assume that $C_{1.5}=\{z\in\mathbb{C}: |z|=1.5\}$ and $g\mid_{C_{1.5}}$ is the rigid rotation.

\medskip
\noindent {\bf Claim 1}.  For each $n\in\mathbb{Z}$, $L_n$ is also a transversal for $\mathcal{P}$ and $L_m\cap L_n=\emptyset$ for any $m\neq n$.
 \begin{proof}[Proof of Claim 1]
 Since $g$ is a homeomorphism and each $C\in\mathcal{P}$ is $g$-invariant, we conclude that $L_n$ is also a transversal.

 Now suppose that there is some $x\in L_{m}\cap L_n$. Let $C\in\mathcal{P}$ be such that $x\in C$. Then we have $g^{-m}x, g^{-n}x\in C\cap L_{0}$. Thus $g^{-m}x=g^{-n}x$, since $|L_0\cap C|=1$. But this contradicts the minimality of $g\mid_{C}$. Thus $L_{m}\cap L_n=\emptyset$.
 \end{proof}

It is easy to see that $\mathcal{P}\mid_{[L_m, L_n]}$ is a partition of $[L_m, L_n]$ for each $m\neq n\in\mathbb{Z}$, where $[L_m, L_n]$ is any region between $L_m$ and $L_n$ that is the closure of any one of the components of $\mathbb{A}\setminus(L_{m}\cup L_{n})$. Precisely, for each $\alpha\in[1,2]$, $[L_m, L_n]\cap C_{\alpha}$ is a curve joining $x_{m}^{\alpha}$ and $x_{n}^{\alpha}$. We denote this curve by $C_{m,n}^{\alpha}$.
\medskip

\noindent {\bf Claim 2}. For each $m\neq n\in\mathbb{Z}$, $\{C_{m,n}^{\alpha}\}_{\alpha\in[1,2]}$ is a continuous decomposition of $[L_{m}, L_{n}]$.

\begin{proof}[Proof of Claim 2]
It suffices to show that $ \{C_{m,n}^{\alpha}\}_{\alpha\in[1,2]}$ is closed in $2^{\mathbb{A}}$. For this, suppose that $C_{m,n}^{\alpha_{i}}\rightarrow K$ in $2^{\mathbb{A}}$ as $i\rightarrow \infty$. By passing to some subsequence, we may assume that $C_{\alpha_{i}}\overset{\prec}{\rightarrow} C_{\alpha}$. We have shown in previously that $C_{\alpha_i}\rightarrow C_{\alpha}$ under the Hausdorff topology. Thus $K\subset C_{\alpha}$. What remains to show is $K=C_{m,n}^{\alpha}$. Since $[L_m, L_n]$ is closed in $\mathbb{A}$, it is clear that $K\subset C_{m,n}^{\alpha}$. On the other hand, it follows from the continuity of $L_{m}$ and $L_n$ that $x_{m}^{\alpha_i}\rightarrow x_{m}^{\alpha}$ and $x_{n}^{\alpha_i}\rightarrow x_{n}^{\alpha}$. Finally, note that $K$ is connected. Thus $K=C_{m,n}^{\alpha}$ and we complete the proof.
\end{proof}

\medskip

\noindent {\bf Claim 3}. If $x_{n_{i}}^{\alpha}\rightarrow x$, then for each $\beta\in [1,2]$, $x_{n_i}^{\beta}\rightarrow y_{\beta}$ for some $y_{\beta}\in C_{\beta}$.
 \begin{proof}[Proof of Claim 3]
 Fix a $\beta\in[1,2]$. To the contrary, assume that there are subsequence $\{k_i\}$ and $\{l_i\}$ of $\{n_i\}$ such that
 \[ x_{k_i}^{\beta}\rightarrow y', \ \ x_{l_i}^{\beta}\rightarrow y'', \ \ \text{ with } \ \ y'\neq y''.\]
 There there are $a,b,c,d\in\mathbb{Z}$ such that there is a  component $(L_{a}, L_{b})$ of $\mathbb{A}\setminus(L_a\cup L_b)$ and a  component $(L_{c}, L_{d})$ of $\mathbb{A}\setminus(L_c\cup L_d)$ with $y'\in (L_a, L_b), y''\in(L_c, L_d)$ and $[L_a, L_b]\cap [L_c, L_d]=\emptyset$, where $[L_a,L_b]=\overline{(L_a,L_b)}$ and $[L_c,L_d]=\overline{(L_c,L_d)}$.
 Then both $[L_a, L_b]\cap\{ x_{n_i}^{\alpha}\}$ and $[L_c, L_d]\cap\{ x_{n_i}^{\alpha}\}$ are infinite. But this contradicts the convergence of $\{x_{n_i}^{\alpha}\}$.
 \end{proof}

Now for each $x\in \mathbb{A}$, let $C_{\alpha}$ be such that $x\in C_{\alpha}$ and take $x_{n_i}^{\alpha}\rightarrow x$. Then let
\[ L_{x}:=\{ y: x_{n_i}^{\beta}\rightarrow y, \beta\in[1,2]\}.\]

\noindent {\bf Claim 4}. $L_x$ is a transversal.
 \begin{proof}[Proof of Claim 4]
 Clearly,  $|L_{x}\cap C_{\beta}|=1$ for any $\beta\in[1,2]$. It remains to show that $L_x$ is an arc.  For this, it suffices to show that the map $[1,2]\rightarrow L_{x}, \alpha\mapsto x_{\alpha}$ is continuous, where $\{x_{\alpha}\}=  L_x\cap C_{\alpha}$.

Fix $\alpha\in (1,2)$ and a neighborhood $U$ of $x_{\alpha}$ in $\mathbb{A}$.  Further, we can take  an open disc $V$ around $x_{\alpha}$ contained in $U$. Take $m,n\in\mathbb{Z}$ such that $C_{m,n}^{\alpha}\subset V$. Then $V$ is a neighborhood of $C_{m,n}^{\alpha}$. By Claim 2,  there are $1<\beta_1<\alpha<\beta_2<2$ such that $C_{m,n}^{\beta}\subset V$ for each $\beta\in[\beta_1,\beta_2]$. In particular, $L_x\cap C_{\beta}\subset V$, for any $\beta\in[\beta_1,\beta_2]$. This shows that $\alpha\mapsto x_{\alpha}$ is continuous.
 \end{proof}

\noindent {\bf Claim 5}. The definition of $L_x$ is independent of the choice of $(n_i)$.
\begin{proof}[Proof of Claim 5]
Suppose that $x_{n_i}^{\alpha}\rightarrow x$ and $x_{m_i}^{\alpha}\rightarrow x$. The we have to show that
\[\lim_{n_i\rightarrow\infty} x_{n_i}^{\beta}=\lim_{m_i\rightarrow \infty} x_{m_i}^{\beta}, \ \ \forall \beta\in[1,2].\]
Let $(k_i)$ be the sequence by putting $(n_i)$ and $(m_i)$ together. Then we have $x_{k_i}^{\alpha}\rightarrow x$. By Claim 3, the sequence $(x_{k_i}^{\beta})$ is convergent for each $\beta$. This implies our Claim.
\end{proof}

 Claim 5 tells us that for each  $y\in L_x$, we have $L_{x}=L_{y}$. Thus $L_{x}, x\in C_{\alpha}$ forms a decomposition of  $\mathbb{A}$. Actually, similar to the proof of the continuity of $\mathcal{P}$, we can also show that $L_{x}, x\in C_{\alpha}$  is a continuous decomposition.
\medskip

\noindent {\bf Claim 6}. For each $\alpha\in[1,2]$, $(x_{n_i}^{\alpha})$ converges in $C_{\alpha}$ if and only if $(n_i\alpha)$ converges in $\mathbb{S}^{1}$.

\begin{proof}[Proof of Claim 6]
This is followed from Claim 3 and our assumption that $C_{1.5}=\{z\in\mathbb{C}: |z|=1.5\}$ and $g\mid_{C_{1.5}}$ is the rigid rotation.
\end{proof}

Now we are ready to show that $g$ can be linearized. Let $\theta$ be the rotation number of $g$. Take a homeomorphism $\psi: L_0\rightarrow P_0:=[1,2]\times \{0\}$.

We define the conjugacy $\Psi: \mathbb{A}\rightarrow \mathbb{A}$ by
\[x\mapsto e^{2\pi i n\theta}\psi(g^{-n}x), \ \ \text{for each }\ \ x\in L_n , n\in\mathbb{Z}, \]
and for $x=\lim x_{n_i}^{\alpha}$, define
\[ \Psi(x)=\lim \Psi(x_{n_i}^{\alpha}).\]

Clearly, we have the following claim.
\medskip

\noindent {\bf Claim 7}. $\Psi g= R_{\theta}\Psi$.
\medskip

In addition, it is also clear that $\Psi$ is one-one. Thus it remains to show the continuity of $\Psi$. For this, suppose that $x_{k}\rightarrow x$ and we assume that $x_{k}\in C_{\alpha_{k}}$. Then we have $\alpha_{k}\rightarrow \alpha$ and $x\in C_{\alpha}$.  We claim that $\lim \Psi(x_k)=\Psi(x)$. WLOG, we may assume that $\lim\Psi(x_k)=y$.  For any neighborhood $V$ of $\Psi(x)$, there are some $m,n\in\mathbb{Z}$ and $\beta_1\prec\alpha\prec \beta_2$ such that $\Psi([\beta_1,\beta_2, L_{m}, L_n])\subset V$, where $[\beta_1,\beta_2, L_{m}, L_n]$ is the  closure of the component of $\mathbb{A}\setminus (\beta_1\cup\beta_2\cup L_{m}\cup L_n)$ containing $x$. Since $x_k\rightarrow x$, $x_{k}\in [\beta_1,\beta_2, L_{m}, L_n]$ for all sufficiently large $k$. This shows that $\Psi(x_k)\rightarrow \Psi(x)$. Thus we have shown that $g$ is conjugate to the rigid rotation by $\Psi$.

\medskip
This completes the proof of Theorem \ref{main2}.

\section{Proof of Theorem \ref{example}}

\subsection{Some notions} Let $\mathbb{A}=\{z\in\mathbb{C}:1\leq |z|\leq 2\}$ be the standard closed annulus on the complex plane. Let $\widetilde{\mathbb{A}}=\{x+yi\in\mathbb{C}:1\leq y\leq 2\}$ be the universal cover of $\mathbb{A}$ and $\pi: \widetilde{\mathbb{A}}\rightarrow \mathbb{A}: x+yi\mapsto ye^{2\pi i x}$ be the covering map. Further, for each $r\in[1,2]$, we denote
\[C_{r}=\{z\in\mathbb{C}: |z|=r\}.\]

Now for each point $z\in \mathbb{A}$, we use the polar coordinate $(\theta, r)$ to represent it  with $\theta\in[0,2\pi)$ and $1\leq r\leq 2$. For each $\beta\in \mathbb{R}$, we define the rigid rotation on $\mathbb{A}$ with angle $2\pi\beta$ by
\[ R_{\beta}: \mathbb{A}\rightarrow\mathbb{A}, z\mapsto ze^{2\pi i\beta}.\]

Let $\gamma: [0,1]\rightarrow \mathbb{A}$ be an essential simple closed curve in $\mathbb{A}$. We define the {\it maximal folding angle} of $\gamma$ to be
\[MFA(\gamma)=\max\left\{|\theta(\gamma(t_1))-\theta(\gamma(t_2))|: \exists 0\leq t_1<t_2<t_3<1 \text{ with }\gamma(t_1)=\gamma(t_3)\right\}.\]

\subsection{Construction of a periodic folding rotation}
We fix $p,q\in\mathbb{N}$ with $(p,q)=1$ and $q> 6$.  We are going to construct a homeomorphism $H$ on $\mathbb{A}$ such that
\begin{enumerate}
\item[(1)] $H R_{p/q}H^{-1}=R_{p/q}$;
\item[(2)] For any $ z,z'\in\mathbb{A}$   with $ |z|=|z'|$, one has
\[ |\theta(H(z))-\theta(H(z'))|\leq 5q  |\theta(z)-\theta(z')|;\]
\item[(3)] the maximal folding angle of $H(C_{3/2})$ is greater than $\frac{\pi}{3}$.
\end{enumerate}

Let $z_0=i, z_1=1+2i, z_2=1+\frac{1}{q}+2i, z_3=\frac{1}{q}+i$, and let
\[ w_1=\frac{1}{2}+\frac{3}{2}i, \ w_2=\frac{3}{4}+\frac{1}{2q}+\frac{7}{4}i, \ w_3=\frac{1}{2}+\frac{1}{q}+\frac{3}{2}i.\]
Define $L_0$ be the union of the segment between $w_1,w_2$ and the  segment between $w_2,w_3$. Then define $ \tilde{L}=\bigcup_{n\in\mathbb{Z}} (L_0+\frac{n}{q})$ and let $\widetilde{\Gamma}=\pi(\tilde{L})$ be the simple closed curve in $\mathbb{A}$.
\medskip

\noindent{\bf Claim 1}. $MFA(\widetilde{\Gamma})> \frac{\pi}{3}$.
\begin{proof}
Note that
\[ Re(w_2)-Re(w_1)=\frac{1}{4}+\frac{1}{2q} \text{ and } Re(w_2)-Re(w_3)=\frac{1}{4}-\frac{1}{2q}.\]
Then we have
\[MFA(\widetilde{\Gamma}) \geq |\theta(\pi(w_2))-\theta(\pi(w_3))| =2\pi\left(\frac{1}{4}-\frac{1}{2q}\right)> \frac{\pi}{3},\]
since $q> 6$.
\end{proof}

Now we are going to construct the desired homeomorphism $H$ on $\mathbb{A}$.  Let $B_0$ be the rectangle whose vertices are $i, 2i, \frac{1}{q}+2i, \frac{1}{q}+i$ and $B_0'$ be the rectangle whose vertices are $z_0, z_1, z_2,z_3$.
First, we define  a homeomorphism $\widetilde{H}$ on $\widetilde{\mathbb{A}}$ such that
\begin{enumerate}
\item[(i)] $\widetilde{H}$ maps $B_0$ onto $B_0'$;
\item[(ii)] $\widetilde{H}$ maps affinely the segment between $\frac{3i}{2} $ and $\frac{1}{q}+\frac{3i}{2}$ onto $L_0$;
\item[(iii)] $\widetilde{H}$ maps affinely each horizontal segment in $B_0$ to a piecewise linear segment in $B_0'$ such that the images of the end points have the same imaginary part (see FIGURE 1);
\item[(iv)] $\widetilde{H}$ is extended to $\widetilde{\mathbb{A}}$ by the translation $z\mapsto z+\frac{1}{q}$, precisely,
\[ \widetilde{H}(z+\frac{n}{q})=\widetilde{H}(z)+\frac{n}{q}, \forall z\in B_0, n\in\mathbb{Z}.\]
\end{enumerate}

Then $\widetilde{H}$ commutes with the integral translations and thus it factors  through a homeomorphism $H$ on $\mathbb{A}$. It is clear that such $H$ satisfies our requirements (1) and (3). It remains to show that $H$ obeys (2). Due the periodicity of $\widetilde{H}$, we may assume that $z,z'\in B_0$ with $|z|=|z'|$.  But then it follows from our construction that
\[ |\theta(H(z))-\theta(H(z'))|\leq \frac{\sqrt{1+(1+1/q)^2}}{1/2q} |\theta(z)-\theta(z')|\leq 5q |\theta(z)-\theta(z')|,\]
since $q>6$.

 It follows from (1) and (2) that we have

\noindent{\bf Claim 2}. For any $z\in\mathbb{A}$,
\[|\theta(HR_{\alpha}H^{-1}(z))-\theta(R_{p/q}(z))|\leq 5q\pi \left|\alpha-\frac{p}{q}\right|.\]

\subsection{Final construction}
We fix an irrational number $\alpha$. It follows from Dirichlet's approximation theorem that there is a sequence $(\frac{p_n}{q_n})$ of rational number with $(p_n,q_n)=1$ and $q_n\rightarrow \infty$ such that
\[ \left|\alpha-\frac{p_n}{q_n}\right|<\frac{1}{q_n^2}, \ \ \forall n\in\mathbb{N}.\]
WLOG, we may further assume  $q_n>6$ for each $n\in\mathbb{N}$.

Now for each $\frac{p_n}{q_n}$, the construction above gives us a homeomorphism $H_n$ on $\mathbb{A}$ satisfying (1), (2), (3).

For each $n\in \mathbb{N}$, let
\[ A_{n}=\left\{z\in\mathbb{C}: 2-\frac{1}{2n}\leq|z|\leq 2-\frac{1}{2n+1}\right\}.\]
Now we define the desired homeomorphism $g$ on $\mathbb{A}$ as follows.
\begin{enumerate}
\item[(a)] The restriction of  $g$ on $A_n$ is conjugate to $H_nR_{\alpha}H_n^{-1}$  by squeezing along the radical direction.
\item[(b)] The definition of $g$ on $\mathbb{A}\setminus \bigcup_{n=1}^{\infty}A_n$ is the rotation $R_{\alpha}$.
\end{enumerate}
It follows from the construction that $g$ is continuous on $\mathbb{A}\setminus C_2$. Thus it suffices to show that  $g$ is continuous on $C_2=\{z\in\mathbb{C}:|z|=2\}$. For this, let $(z_k)$ be a sequence in $\mathbb{A}$ with $z_k\rightarrow z\in C_2$. Since the restriction of $g$ outside $\bigcup_{n=1}^{\infty}A_n$ is the standard rotation $R_{\alpha}$. WLOG, we may assume that $\{z_k\}\subset\bigcup_{n=1}^{\infty}A_n$. Further, to show that $g(z_k)\rightarrow g(z)=ze^{2\pi i\alpha}$ it suffices to show that
\[\theta(g(z_k))-\theta(z_k)\rightarrow 2\pi\alpha.\]
We assume $z_k\in A_{n_{k}}$. We denote the conjugation from $g|_{A_n}$ to $H_nR_{\alpha}H_{n}^{-1}$ by $\phi_n$, which is a squeezing transformation along the radical direction. Then we have
\begin{align*}
 \theta(g(z_k))-\theta(z_k)&=\theta(H_{n_k}R_{\alpha}H_{n_k}^{-1}(\phi_{n_k}(z_k)))-\theta((\phi_{n_k}(z_k))\\
  &=\theta(H_{n_k}R_{\alpha}H_{n_k}^{-1}((\phi_{n_k}(z_k)))-\theta(R_{p_{n_k}/q_{n_k}}((\phi_{n_k}(z_k)))+2\pi(p_{n_k}/q_{n_k})\\
  &\leq 5q_{n_k}\pi\left|\alpha-p_{n_k}/q_{n_k} \right| +2\pi(p_{n_k}/q_{n_k})\\
  &\leq \frac{5\pi}{q_{n_k}}+2\pi(p_{n_k}/q_{n_k}).
\end{align*}
Thus $\theta(g(z_k))-\theta(z_k)\rightarrow 2\pi\alpha$.

It remains to show that $g$ cannot be linearizable. For this, it suffices to show that there is no a transversal. But this follows from the property (3) of $H$ that there is a folding circle in each $A_n$ whose maximal folding angle is greater than $\frac{\pi}{3}$.

 \begin{figure}[http]
 \includegraphics[width=10cm]{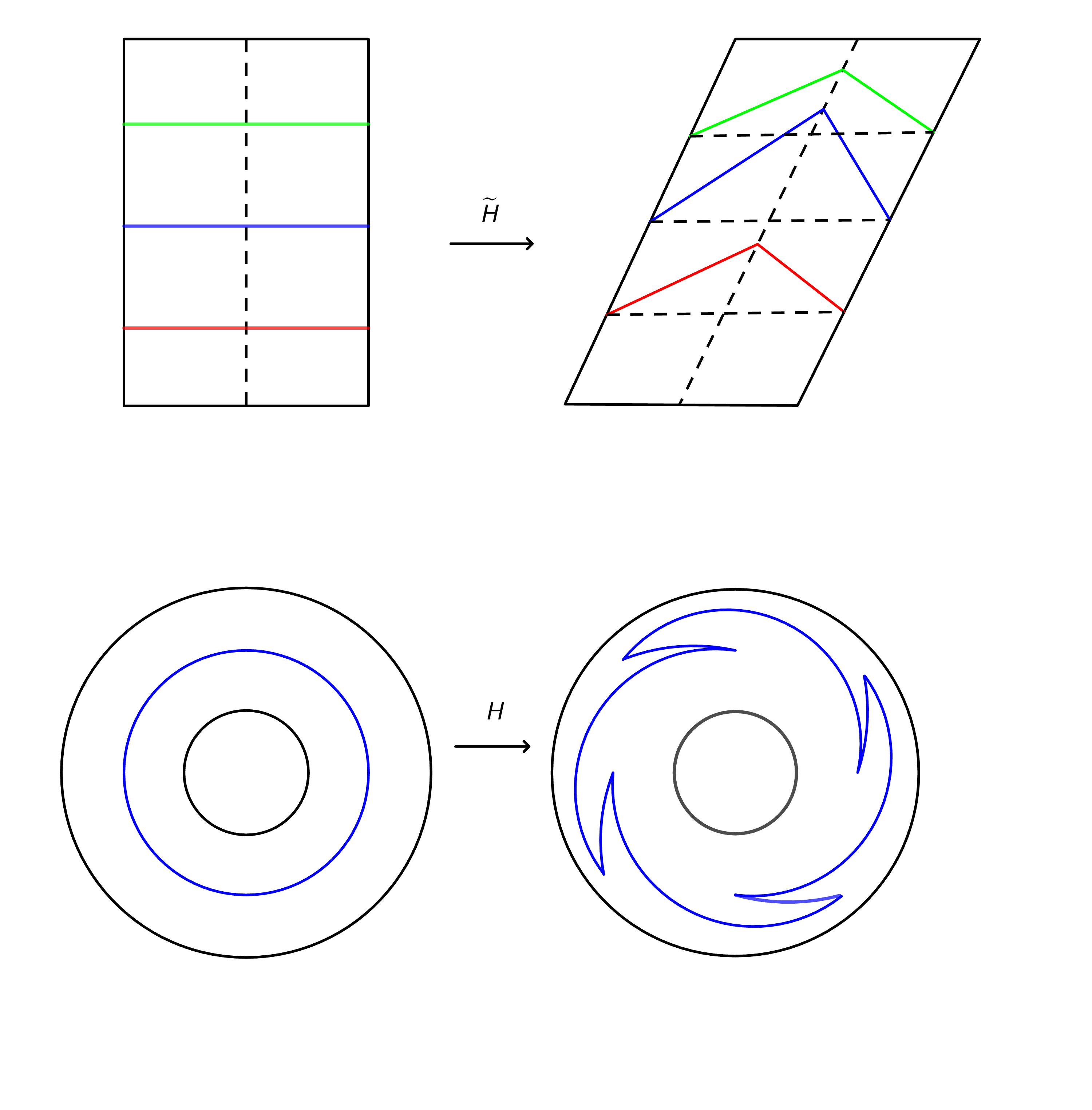}
 \caption{ }
 \end{figure}


\end{document}